\documentclass{conm-p-l}
\usepackage{amsfonts,amstext,amsmath}
\begin{document}
\newcommand{\cal}{\mathfrak}
\newtheorem{Th}{Theorem}
\newtheorem{Def}{Definition}
\newtheorem{Lemm}{Lemma}
\newtheorem{Cor}{Corollary}
\newtheorem{Rem}{Remark}
\newtheorem{Exam}{Example}
\newtheorem{Prop}{Proposition}

\title[Miraculous Cancellation]{\bf Miraculous Cancellation and Pick's Theorem.}
\author{K.~E.~Feldman
}
\address{DPMMS, University of Cambridge,
Wilberforce Road, Cambridge, UK, CB3 0WB}
\email{k.feldman@dpmms.cam.ac.uk}
\thanks{The research
was supported by EPSRC grant GR/S92137/01.}
\begin{abstract}
We show that the Cappell--Shaneson version of Pick's theorem for
simple lattice polytopes is a consequence of a general relation between
characteristic numbers of virtual submanifolds dual to the
characteristic classes of a stably almost complex manifold. This
relation is analogous to the miraculous cancellation formula of
Alvarez-Gaume and Witten, and is imposed by the action of the
Landweber--Novikov algebra in the complex cobordism ring of a
point.
\end{abstract}
\maketitle


\section*{Introduction}
Studying gravitational anomalies Alvarez-Gaume and
Witten~\cite{AlWi} discovered a remarkable cancellation which
they called the miraculous cancellation formula. The formula is a
consequence of the following relation between Pontryagin
characteristic numbers of 12-dimensional oriented manifolds:
\begin{equation}
\label{MirCan}
L(M)=8A(M,T)-32\hat A(M),
\end{equation}
where $L(M)$ is the $L$-genus (the signature) of $M$ corresponding
to the power series $L(x)=x/\tanh(x)$, $\hat A(M)$ is the $\hat
A$-genus of $M$ corresponding to the power series
$A(x)=x/2\sinh(x/2)$, and $A(M,T)$ is the twisted $\hat A$-genus
defined by cohomology Kronecker product
$$
\langle \prod_jA(x_j)\left(\sum_j
e^{x_j}+e^{-x_j}\right),[M]\rangle
$$
($\pm x_j$ are the formal Chern roots of $\mathbb C\otimes \tau(M)$ in
ordinary cohomology and $[M]\in H^{\dim M}(M,\mathbb Q)$ is the
fundamental class).

The miraculous cancellation formula was generalized by Kefeng
Liu~\cite{Liu} to higher dimensions with the use of modular
invariance of certain elliptic operators on the loop space, it was
strengthen further in~\cite{FeZh}.

Motivated by this development Buchstaber and Veselov suggested
in~\cite{BV} that a general formula exists which relates
characteristic numbers of virtual submanifolds Poincare dual to
cobordism characteristic classes of a given manifold. In
particular, they found that the signatures of manifolds Poincare
dual to the tangent Pontryagin classes with values in complex
cobordism of an oriented manifold $M^{4n}$ are related as
$$
L(M^{4n})+\sum^n_{k=1}(-1)^k
L\left(\left[P^U_k(\tau(M^{4n}))\right]\right]=0 \left(\rm{mod}
2^{\alpha(n)}\right),
$$
where $\alpha(n)$ is the number of ones in the binary expansion of
the number $n$. The general formula which expresses all relations
between Chern characteristic numbers of virtual submanifolds was
obtained in~\cite{Feldman}.

Various applications of miraculous cancellation type formulae to
the questions of divisibility of topological invariants were found
in~\cite{FaOs, Feldman, FeZh, LZ}. In this paper we show that
Cappell--Shaneson version~\cite{CappellShaneson} of Pick's theorem is, in fact, 
also a consequence
of a cancellation formula for the Chern characteristic numbers of
virtual submanifolds.

 Classical Pick's theorem~\cite{Pick} states that the area of a convex polygon
$P$ with all vertexes in the standard two-dimensional unit lattice
$\mathbb Z^2\subset \mathbb R^2$ can be expressed in terms of the
number of lattice points inside the polygon as
$$
Area(P)=Int(P)+\frac{Bd(P)}{2}-1,
$$
where $Int(P)$ is the number of lattice points strictly inside $P$
and $Bd(P)$ is the number of lattice points on the boundary of
$P$. Numerous generalizations (see~\cite{BriVer, Gu1, KSW, KP,
Morelli, Pommersheim} and references there) of this statement to
higher dimensions are based on the theory of toric
varieties~\cite{Danilov}. Toric varieties are algebraic varieties
naturally associated with a simple lattice polytope in $\mathbb
R^n$. The twisted Todd genus of a toric variety is the number of
lattice points in the corresponding lattice polytope. Lack of
elementary description of the Todd class results in complexity of
combinatorial formulae for the number of lattice
points~\cite{Morelli, Pommersheim}.

 One may look at the classical Pick's theorem from a
slightly different angle.  We may argue that in Pick's theorem
lattice points should be counted with different weights depending
on the dimension of faces they belong to. Then for a lattice
$n$-gon $P$ Pick's theorem takes the form
$$
Area(P)=Int(P)+\frac{Edg(P)}{2}+\frac{Vert(P)}{4}-\frac{4-n}{4},
$$
where $Edg(P)$ is the number of lattice points inside each edge of
$P$ and $Vert(P)=n$ is the number of vertexes of $P$. Cappell and
Shaneson observed that for a general simple lattice polytope
$P\subset\mathbb R^n$ the sum
$$
\sharp({\rm rel.int}(P)\bigcap \mathbb
Z^n)+\sum^n_{k=1}\left(\frac{1}{2}\right)^k \sum_{{F\subset
P}\atop {\dim F=(n-k)}}\sharp\left({\rm rel.int}(F)\bigcap \mathbb
Z^n\right)
$$
can be expressed in terms of the polytope algebra over the ring of
infinite order constant coefficient differential operators on
$\mathbb R^n$. In this paper we prove that for any stably almost
complex manifold $M^{2n}$ with some fixed cohomology class $w\in
H^2(M^{2n},\mathbb Z)$
$$
\langle exp(w){\rm Td}(M),[M]\rangle
+\sum^n_{k=1}\left(-\frac{1}{2}\right)^k \langle exp(w_k){\rm
Td}([c_k(\tau(M))]),[c_k(\tau(M))]\rangle=
$$
$$
=\langle exp(w)\prod^m_{i=1}
\frac{x_i/2}{\tanh(x_i/2)},[M]\rangle,
$$
where $[c_k(\tau(M))]$ are manifolds Poincare dual to the
cobordism Chern characteristic classes of $\tau(M)$, ${\rm Td}$ is
the Todd class, $x_i$ are the tangent Chern roots of $M$ and
$w_k\in H^2([c_k(\tau(M))],\mathbb Z)$ are the pullbacks of $w$
under the canonical embeddings $[c_k(\tau(M))]\subset M\times
\mathbb R^N$. In the case of toric varieties this formula displays
a duality between cobordism invariants of a smooth toric variety
and its canonical submanifolds defined by faces of the
corresponding polytope. Twisted Todd genera of these submanifolds
are exactly the numbers of lattice points inside the closures of
the corresponding faces and the Pick's type theorem reads as
\begin{equation}
\label{CaSh}  \sharp(P\bigcap \mathbb
Z^n)+\sum^n_{k=1}\left(-\frac{1}{2}\right)^k \sum_{{F\subset
P}\atop {\dim F=(n-k)}}\sharp\left(F\bigcap \mathbb Z^n\right)
=\int_{M_P} exp(w_P)\prod^m_{j=1}\frac{x_i/2}{\tanh(x_i/2)},
\end{equation}
where $w_P$ is the Kahler class of a toric variety $M_P$
corresponding to the simple lattice polytope $P$, and $x_i$, $i=1,\dots,
m$, are Chern classes of the canonical line bundles corresponding
to the codimension one faces of $P$. The right hand side
of~(\ref{CaSh}) is a weighted analogue of the twisted signature of
$M_P$.

In lower dimensions the twisted signature can be developed
explicitly. We use it for evaluation of a Pick's type theorem for
a three-dimensional example. As an intermediate step in order to
compute the twisted signature we obtain a general expression for
monomial Chern characteristic numbers of a smooth toric variety of
any dimension in terms of primitive vectors which define
codimension one faces of the polytope.

The paper is organized as follows. The first three sections give a
brief account of the facts used for the proof of the miraculous
cancellation type formula. All information in these sections is
well know to experts and is given mainly to set up the notations.
The cancellation formula itself is proved in Section 4. In Section 5 we show how
the cancelation formula of Alvarez-Gaume and Witten can be deduced
within the same cobordism technique.
In Section
6 we state basic facts about toric manifolds and give an
expression for the Gysin map in terms of the fixed point data.
Finally, in Section 7 we derive Cappell--Shaneson type formula and
evaluate Pick's theorem for some low--dimensional examples.


\section{Bordism Ring of Manifolds with Line Bundles}

 Let us consider a set of all pairs consisting of a compact stably
almost complex manifold  without boundary and a complex line
bundle over it. We introduce an equivalence relation on this set
by means of a standard construction from cobordism theory. We say
that two pairs $(M^n_1,\eta_1)$ and $(M^n_2,\eta_2)$ are
equivalent if there is a compact stably almost complex manifold
$W$ of dimension $n+1$ with a complex line bundle $\eta$ over it
such that
$$
\partial W = M^n_1\bigcup M^n_2,\qquad \eta|_{M^n_1}\cong \eta_1,\qquad
\eta|_{M^n_2}\cong\eta_2,
$$
and the restrictions of the almost complex structure in the normal
bundle $\nu(W)$ of $W$ on $M^n_i$, $i=1,2$, coincide with
$\nu(M^n_1)$ and $\overline{[1]}_{\mathbb C}\oplus \nu(M^n_2)$
respectively ($[1]_{\mathbb C}$ is the trivial complex line
bundle). Where it is necessary we write $(M,\nu(M),\eta)$ instead
of $(M,\eta)$ to underline the fixed (in the normal bundle) stably
almost complex structure on $M$.

 The set of pairs $(M,\eta)$ factorized by the above equivalence relation,
is a graded Abelian group with an addition given by a disjoint
union of manifolds. The inverse to $(M,\nu(M),\eta)$ is
$(M,\overline{[1]}_{\mathbb C}\oplus\nu(M),\eta)$. Zero is the
empty set. We refer to this group as $\Omega^{LU}_*$.

 We introduce a multiplication in $\Omega^{LU}_*$ by the formula given
on representatives of bordism classes:
$$
[(M_1,\eta_1)]\cdot [(M_2,\eta_2)]=[(M_1\times
M_2,\eta_1\otimes\eta_2)].
$$
Obviously, $\Omega^{LU}_*$ equipped with such an operation becomes
a commutative ring with a unit $(pt,[1]_{\mathbb C})$. We can
easily calculate the ring $\Omega^{LU}_*$.
\begin{Prop}
\label{U(CP)} The graded Abelian group $\Omega^{LU}_{**}$ is
isomorphic to the complex bordism group of the infinite complex
projective space:
$$
\Omega^{LU}_*\cong U_*(CP(\infty),\emptyset).
$$
Therefore, the ring $\Omega^{LU}_*$ is the polynomial ring $\Omega^U_*[[t]]$, where $\Omega^U_*$
is the complex bordism ring (complex cobordism group of a point).
\end{Prop}
 For any pair of topological spaces $(X,A)$ with basepoints we define the
corresponding $\Omega^{LU }_*$-bordism group $LU_*(X,A)$ which
consists
 of all triples $(M,\eta,f)$, where $M$ is a compact stably almost complex
manifold with boundary, $\eta$ is a complex line bundle over $M$
and $f$ is a continuous map:
$$
f:(M,\partial M)\to (X,A),
$$
factorized by the bordism relation:
$$
(M_1,\eta_1,f_1)\sim (M_2,\eta_2,f_2),
$$
if and only if there is a triple $(W,\eta,F)$, such that
$$
\partial W = M_1\bigcup M_2\bigcup V,\qquad
F:(W,V)\to (X,A),
$$
$$
F|_{M_i}=f_,\qquad \nu|_{M_1}=\nu_1,\qquad
\nu|_{M_2}=\overline{[1]}_{\mathbb C}\oplus \nu_1.
$$

Denote the Thom space of the tautological complex vector bundle
$\eta_n$ over $BU(n)$ by $MU(n)$. Let $E_{2k}$, $k=0,1,...$, be
the wedges $CP(\infty)_+\wedge MU(k)$, where $X_+$ is a disjoint
union of a topological space $X$ and a point $\{pt\}$, we will
also use an identification $X_+=(X,\emptyset)$. Set
$E_{2k+1}=\Sigma E_{2k}$. Define a sequence of maps:
$$
f_k:\Sigma E_k\to E_{k+1}
$$
by the formulae
$$
f_{2k}=id,
$$
$$ 
f_{2k+1}:\Sigma
E_{2k+1}=\Sigma^2E_{2k}=CP(\infty)_+\wedge \Sigma^2
MU(k)\stackrel{id\wedge j_k}\longrightarrow  CP(\infty)_+\wedge
MU(k+1),
$$
where $j_k:\Sigma^2 MU(k)\to MU(k+1)$ is the standard embedding.

 There is a sequence of canonical maps  $i_n:S^n\to E_n$ which is inherited
from the canonical embedding $S^{2n}\subset MU(n)$. We also define
a graded paring for any two spaces $E_k$ and $E_m$:
$$
\mu_{km}:E_k\wedge E_m\to E_{k+m},
$$
by the formula given on even spaces $E_n$:
$$
\mu_{2k,2m}:E_{2k}\wedge E_{2m}= \left(CP(\infty)\times
CP(\infty)\right)_+\wedge MU(k)\wedge MU(m)
\stackrel{\otimes\wedge i}\longrightarrow 
$$
$$
\to
CP(\infty)_+\wedge
MU(k+m),
$$
where $\otimes:CP(\infty)\times CP(\infty)\to CP(\infty)$ is the
classifying map for the tensor product of Hopf's line bundles over
$CP(\infty)$, $i:MU(k)\wedge MU(m)\to MU(k+m)$ is the canonical
embedding.
\begin{Prop}
\label{coh.th} The set of spaces $\{E_k\}^{\infty}_{k=0}$ and maps
$\{f_k,j_l,\mu_{mn}\}$ forms a multiplicative spectrum, and, thus,
defines a (co)homology theory:
$$
h_s(X,A)=\lim_{t\to \infty}\pi_{s+t}(X/A\wedge E_t),
$$
$$
h^s(X,A)=\lim_{t\to \infty}[\Sigma^t(X/A),E_{s+t}]
$$
The corresponding cohomology theory is multiplicative.
\end{Prop}
By the standard Thom transversality arguments we deduce
\begin{Prop}
\label{Cob} The bordism theory $LU_*(\cdot)$ coincides with the
homology theory $h_*(\cdot)$.
\end{Prop}
We denote the corresponding cobordism theory by $LU^*(\cdot)\cong
h^*(\cdot)$.


\section{Poincare Duality and the Euler Class}

Let us consider a compact stably almost complex manifold $M^n$
without boundary. Assume that this manifold is equipped with a
complex line bundle $\zeta$. Any element $\alpha\in LU^k(M^n)$ can
be realised as a compact stably almost complex submanifold $L$,
$\dim L=n-k$, in $M^n\times {\mathbb R}^N$ ($N$ is large) without
boundary, equipped with a complex line bundle $\zeta_L$ and with a
fixed stable complex structure $\nu$ in the normal bundle of the
embedding $L\subset M^n\times {\mathbb R}^N$ satisfying 
$\nu(L)\cong_{\mathbb C}\nu\oplus \nu(M\times {\mathbb R}^N)|_L$. Denote the embedding
$L\subset M^n\times {\mathbb R}^N$ by $i$ and the projection
$M^n\times {\mathbb R}^N\to M^n$ by $p$. Define a map:
$$
D:LU^*(M^n)\to LU_{n-*}(M^n)
$$
by the formula
$$
D[(L,\nu,\zeta_L)]=[(L,\bar{\zeta_L}\otimes i^*p^*\zeta)]
$$
\begin{Prop}
\label{Dual} For any compact stably almost complex manifold $M^n$
equipped with a line bundle $\zeta$ the map
$$
D=D(M^n;\zeta):LU^*(M^n)\to LU_{n-*}(M^n)
$$
is an isomorphism.
\end{Prop}
\begin{Cor}
\label{FundCl} For any compact closed stably almost complex
manifold $M^n$ equipped with a complex line bundle $\eta$, the
fundamental class of $M^n$, i.e. $D$-dual to $1\in LU^0(M^n)$, is
$$
[M^n,\eta]=[(M^n,id,\eta)]\in LU_n(M^n).
$$
\end{Cor}
Let us consider a complex $n$-dimensional vector bundle $\xi$ over
a base space $X$. Assume that $X$ is equipped with a complex line
bundle $\zeta$.
\begin{Def}
\label{ThomCl} The canonical Thom class of $(\xi,\zeta)$ is an
element $u(\xi,\zeta)\in LU^{2n}(T\xi)$ represented by the map
$$
\alpha:T\xi\stackrel{p\wedge id}\longrightarrow X_+\wedge
T\xi\stackrel{f_{\bar\zeta}\wedge g_{\xi}}\longrightarrow
CP(\infty)_+\wedge MU(n),
$$
where $p(x,v_x)=x$ is the projection outside infinity,
$f_{\bar\zeta}:X\to CP(\infty)$ is the map classifying
$\bar\zeta$, $g_{\xi}:T\xi\to MU(n)$ is the map classifying $\xi$.
The homotopy class of $\alpha$ lies in $[T\xi,CP(\infty)_+\wedge MU(n)]$
and, thus, it defines an element 
$$
[\alpha]\in  \lim_{k\to
\infty}[\Sigma^{2k}T\xi,CP(\infty)_+\wedge MU(n+k)]=LU^{2n}(T\xi);
$$

\end{Def}
\begin{Prop}
For a pair $(\xi,\zeta)$ of two complex vector bundles over a
finite CW-complex $X$ ($\dim_{\mathbb C}\xi=n$, $\dim_{\mathbb
C}\zeta=1$) the multiplication by $u(\xi,\zeta)$
$$
LU^*(X,\emptyset)\longrightarrow LU^{*+2n}(T\xi),
$$
is an isomorphism.
\end{Prop}

\

Consider a complex $n$-dimensional vector bundle $\xi$ over a
topological space $X$. Assume that $X$ is equipped with a complex
line bundle $\zeta$. Denote the classifying maps for $\xi$ and
$\bar\zeta$ by $f$, $g$ respectively:
$$
f:X\to BU(n),\qquad \xi=f^*\eta_n,\qquad \hat f:T\xi\to MU(n);
$$
$$
g:X\to CP(\infty), \qquad \bar\zeta=g^*\eta_1.
$$
Let $s_0:X\to T\xi$ be the  embedding onto zero section of $\xi$.
\begin{Def}
\label{Euler} The Euler class of the pair $(\xi,\zeta)$ is
$$
\chi(\xi,\zeta)=s^*_0 u(\xi,\zeta);
$$
$$
\chi(\xi,\zeta)=[g\wedge( \hat f\circ s_0)]\in \lim_{t\to \infty}
[\Sigma^{2t}(X_+),CP(\infty)_+\wedge
MU(n+t)]=LU^{2n}(X,\emptyset).
$$
\end{Def}

Once we have defined the Euler class we can continue and introduce
the other characteristic classes for $(\xi,\zeta)$ following
techniques from~\cite{B1,CF}.

Let $\mathbb CG^n_k(\xi)$ be the complex grassmannization of the
complex vector bundle $\xi$. Denote the canonical $k$-dimensional
complex vector bundle over $\mathbb CG^n_k(\xi)$ by $\xi(k)$. Let
$p_k:\mathbb CG^n_k(\xi)\to X$ be the projection on the base space
$X$, and $\tau(p_k):S^N\wedge X_+\to S^N\wedge \mathbb
CG^n_k(\xi)_+$ be the transfer map of Becker--Gottlieb~\cite{BG}
for fibre bundle $(\mathbb CG^n_k(\xi),\mathbb CG^n_k,X,p_k)$.
\begin{Def}
\label{ChCl} Characteristic classes $c^{LU}_{km}(\xi,\zeta)$ with
values in $LU^*(X,\emptyset)$ of $(\xi,\zeta)$ are defined by
$$
c^{LU}_{km}(\xi,\zeta)=\tau(p_k)^*\{\chi(\xi(k),p^*_k\zeta^m)\}
\in LU^{2k} (X,\emptyset),
$$
here $\dim_{\mathbb C}\xi=n$, $\dim_{\mathbb C}\zeta=1$.
\end{Def}


\section{The Chern-Dold character}

For any multiplicative cohomology theory $h^*(\cdot)$ there is a
unique multiplicative transformation
$$
ch_h:h^*(\cdot)\to H^*(\cdot;h^*(pt)\otimes {\mathbb Q}),
$$
such that for a point it is the identical embedding. Tensored with
the rational numbers this transformation is an isomorphism of
cohomology theories.

Consider a compact stably almost complex manifold $M^{2n}$ with a
fixed complex structure in its normal bundle $\nu$, $\dim_{\mathbb
C}\nu=m$. Let   $\zeta$ be a complex line bundle over $M$. Chose
$1\in LU^0(M^{2n},\emptyset)$. This element can be represented by
a triple $(M^{2n},[1]_{\mathbb C},id)$.
\begin{Prop}
\label{chu} Let $p:M^{2n}\to \{pt\}$ be the constant map. In the
notations described above:
$$
[M^{2n},\zeta]= ch_{LU}\circ p^{LU}_!(1)=\langle
T_{LU}(\nu),[M^{2n}]\rangle,
$$
where $p^{LU}_!:LU^*(M^{2n},\emptyset)\to LU^{*-2n}(S^0)$ is
the Gysin map, $T_{LU}$ is the Todd class of the transformation
$ch_{LU}$, $[M^{2n}]$ is the fundamental class of the manifold
$M^{2n}$ in homology $H_*(M^{2n};\Omega^{LU}_*\otimes {\mathbb
Q})$.
\end{Prop}
{\bf Proof.} The first equality follows from
Corollary~\ref{FundCl} and the second equality is the standard
consequence of the Riemann--Roch theorem.

\

For any compact closed stably almost complex manifold $M^{2n}$
equipped with a complex line bundle $\zeta$ we define normal mixed
characteristic numbers:
$$
c_{I}(M,\zeta)=\langle
c_{I}(\nu(M))c^{n-|I|}_1(\zeta);[M^{2n}]\rangle,
$$
where $c_I$ is the Chern class in ordinary cohomology
corresponding to the partition $I=(i_1,...,i_k)$,
$|I|=i_1+...+i_k\le n$. It is a well--known fact~\cite{Stong} that
two compact closed stably almost complex manifolds $M_1$, $M_2$
equipped with complex line bundles $\zeta_1$, $\zeta_2$
respectively represent the same element in $\Omega^{LU}_*$ if and
only if they have the same  mixed characteristic numbers.

Let $S(x)$ be an infinite power series in $x$  defined by the
formula:
$$
S(x)=\prod^{\infty}_{i=1}(1-y_ix)^{-1},
$$
where $y_i$, $i=1,...$, are independent variables. Define a map:
$$
S_t:\Omega^{LU}_*\to \Lambda_y[t],
$$
$$
S_t([(M^{2n},\nu(M^{2n}),\zeta)])=\langle \mathrm
e^{-tc_1(\zeta)}\prod_j S(x_j),[M^{2n}]\rangle
$$
where $\Lambda_y$ is the ring of symmetric polynomials over
$\mathbb Q$ in variables $y_i$, $=1,2,...$, and $x_j$,
$j=1,2,...$, are the Chern roots of the normal bundle
$\nu(M^{2n})$.
\begin{Prop}
\label{Genus} The map $S_t$ is a monomorphism which induces an
isomorphism
$$
S_t:\Omega^{LU}_*\otimes {\mathbb Q}\to \Lambda_y[t].
$$
\end{Prop}
\begin{Cor}
\label{Specialization} For every ring homomorphism
$\phi:\Omega^{LU}_*\to \mathbb Q$ there is such a rational valued
specialization for elementary symmetric functions
$e_i=e_i(y_1,y_2,...)$ and for variable $t$, that $S_t$ for this
specialization coincides with $\phi$.
\end{Cor}
Let $S^*_t$ be the transformation
$$
H^*(\cdot;\Omega^{LU}_*\otimes {\mathbb Q})\to
H^*(\cdot;\Lambda_y[t]),
$$
induced by the coefficient ring homomorphism $S_t$.
\begin{Prop}
The Todd class of the composition $\tau=S^*_t\circ ch_{LU}$ is
$$
T_{\tau}(x,z)=\mathrm e^{-tz}\prod^{\infty}_{i=1}(1-y_ix)^{-1}.
$$
In this expression for a pair $(\xi,\zeta)$ of complex bundles
over the base space $X$ ($\dim_{\mathbb C}\zeta=1$)  the variable
$x$ is fixed for the Chern roots in ordinary cohomology of $\xi$,
and the variable $z$ is fixed for the first Chern class of
$\zeta$.
\end{Prop}

\

Let us calculate $LU^*(CP(n),\emptyset)$. We fix the canonical
line bundle $\eta_1$ over $CP(n)$. Observe also that
$CP(n+1)\approx T\bar\eta_1$. Thus we have

\begin{Lemm}
There is an isomorphism
$$
LU^{*+2}(CP(n+1))\cong LU^*(CP(n),\emptyset),
$$
via multiplication by the Thom class $u(\bar\eta_1,\eta_1)$.
\end{Lemm}

\begin{Cor}
The ring $LU^*(CP(n),\emptyset)$ is a polynomial ring
$\Omega^{LU}_*[u]/\{u^{n+1}=0\}$. It is additively generated by
$[(CP(k),\eta^{n-k}_1,i_k)]$, $i_k:CP(k)\stackrel{\nu_k}
\longrightarrow CP(n)$, $k=1,...,n$. The multiplication is given
by the rule
$$
[(CP(k),\eta^{n-k}_1,i_k)]\cdot[(CP(m),\eta^{n-m}_1,i_m)]
=[(CP(k+m-n),\eta^{2n-k-m}_1,i_{k+m-n})],
$$
\end{Cor}
\begin{Prop}
\label{gener_image} Via the map 
$$
\tau=S^*_t\circ
ch_{LU}:LU^*(CP(\infty),\emptyset)\to
H^*(CP(\infty);\Lambda_y[t])
$$ 
the canonical generator $u\in
LU^*(CP(\infty),\emptyset)$ goes to
$$
x\mathrm e^{-tx}\prod^{\infty}_{i=1}(1-y_ix)^{-1}
$$
where $x$ is the canonical generator of
$H^*(CP(\infty),\Lambda_y[t])$.
\end{Prop}


\section{General Cancellation Formula}

The theory $LU^*$ possesses an algebra of stable operations
similar to that of Landweber--Novikov algebra in complex
cobordism~\cite{Lan,N}. Namely, we can combine them in one
multiplicative operation which on the Thom class
$u(\eta_n,\eta_1)$ of $\left(CP(\infty)_+\right)\wedge MU(n)$
takes the following value:
$$
{\cal S}_{(z)}\left(u(\eta_n,\eta_1)\right)= u(\eta_n,\eta_1)
\left(1+\sum_{\lambda>0} h_{\lambda}(z)m_{\lambda}(\eta_n)\right),
$$
where $h_{\lambda}(z)$ is the complete symmetric function
corresponding to the partition
$\lambda=(\lambda_1\ge\lambda_2\ge\cdots\ge 0)$~\cite{Mac},
$m_{\lambda}(\eta_n)$ is the $LU^*$-Chern class corresponding to
the monomial symmetric function $m_{\lambda}$ and a trivial line
bundle.
 Similar to~\cite{Feldman} we conclude that
\begin{Prop}
The composition $\Theta=S^*_{(y,t)}\circ ch_{LU}\circ {\cal
S}_{(z)}$:
$$
\Theta: LU^*(\cdot)\to H^*\left(\cdot,
\Lambda_y\otimes\Lambda_z\otimes {\mathbb Q}[t]\right)
$$
is a multiplicative transformation of cohomology theories and its
Todd class is
$$
T_{\Theta}(\eta,\zeta)=e^{-tu}\prod^n_{i=1}\left(\prod^{\infty}_{j=1}
(1-y_jx_i)^{-1}\cdot \prod^{\infty}_{k=1}
\left(1-\frac{z_kx_i}{\prod^{\infty}_{l=1}(1-y_lx_i)}\right)^{-1}\right),
$$
where $u$ is the first Chern class of the complex line bundle
$\zeta$, $x_i$, $i=1,...,n$, are the Chern roots in integral
cohomology of the complex vector bundle $\eta$, $n=\dim_{\mathbb
C}\eta$.
\end{Prop}

Applying the Riemann--Roch theorem to the transformation $\Theta$
and to the constant map $p:M\to \{pt\}$ with $LU$-oriented
manifold $(M^{2n},\zeta)$, we deduce the following theorem:

\begin{Th}
\label{formula} Let
$([m_{\lambda}(\nu(M^{2n}))],i^*_{\lambda}\zeta)$ be manifolds
$LU$-dual to the $LU$-Chern classes of the normal bundle
$\nu(M^{2n})$ which correspond to the monomial symmetric functions
$m_{\lambda}$, $\lambda=(\lambda_1\ge\lambda_2\ge\cdots 0)$. Then
$$
\langle T_{\Theta}(\nu(M^{2n}),\zeta),[M^{2n}]\rangle= 
$$
$$
=
\langle
T_{\tau}(\nu(M^{2n}),\zeta),[M^{2n}]\rangle+
\sum_{\lambda>0}h_{\lambda}(z) \langle
T_{\tau}(\nu_{\lambda},i^*_{\lambda}\zeta),[m_{\lambda}(\nu(M^{2n}))]
\rangle,
$$
where $\nu_{\lambda}$ is the normal bundle of $m_{\lambda}(\nu(M^{2n}))$
with the induced stable almost complex structure.
\end{Th}

\

Let $M$ be a stably almost complex manifold and $x_i$,
$i=1,2,...$, be the Chern roots of the stable complex tangent
bundle to $M$. We define the Todd class of $M$ by
$$
{\rm Td}(M)=\prod^{\infty}_{i=1} \frac{x_i}{1-e^{-x_i}}.
$$
This is a symmetric polynomial in $x_i$, $i=1,2,\dots$, which can
be expressed as a series in tangent Chern classes of $M$ in
ordinary cohomology.
\begin{Th}
\label{Mirac} Let $M$ be a stably complex manifold and $w\in
H^2(M,\mathbb Z)$. Then for the manifolds $[c_k(\tau(M))]$ dual to
the complex cobordism Chern classes $c^U_k$ of $\tau(M)$ the
following formula holds
$$
\langle exp(w)\prod^m_{i=1} \frac{x_i/2}{\tanh(x_i/2)},[M]\rangle=
$$
$$
= \langle exp(w){\rm Td}(M),[M]\rangle
+\sum^n_{k=1}\left(-\frac{1}{2}\right)^k \langle exp(w_k){\rm
Td}([c_k(\tau(M))]),[c_k(\tau(M))]\rangle,
$$
where $x_i$ are the tangent Chern roots of $M$ and $w_k\in
H^2([c_k(\tau(M))],\mathbb Z)$ are the pullbacks of $w$ under the
canonical embeddings $[c_k(\tau(M))]\subset M\times \mathbb R^N$.
\end{Th}
{\bf Proof.} Let us substitute the following values for the
variables $t,z_i$ and $y_j$ in the series
$T_{\Theta}(\eta,\zeta)$, where $\eta$ is the normal bundle of
$M$, and $\zeta$ is a line bundle over $M$ with $c_1(\zeta)=w$:
$$
t=1,\quad z_1=1/2,\quad z_2=z_3=...=0,
$$
$$
\prod^{\infty}_{j=1}(1-y_jx)^{-1}=\frac{1-e^{-x}}{x}.
$$
For these values of variables we obtain
$$
T_{\Theta}(\eta,\zeta)=e^{u}\prod^n_{i=1}
\frac{1-e^{-x_i}}{x_i}\frac{2}{1+e^{-x_i}}=e^u\prod^n_{i=1}
\frac{\tanh(x_i/2)}{x_i/2},
$$
where $x_i$ are the Chern roots of $\eta$ and $u=c_1(\zeta)$.
Therefore, for the tangent Chern roots $x'_i$, $i=1,2,\dots$, the
left hand side of the formula in Theorem~\ref{formula} becomes
$$
\langle exp(w)\prod^m_{i=1}
\frac{x'_i/2}{\tanh(x'_i/2)},[M]\rangle.
$$
 Because $c_k(\tau(M))=(-1)^k
h_k(\nu(M))$ for our choice of variables $z_i$, the right hand
side of  Theorem~\ref{formula} is exactly
$$
\langle exp(w){\rm Td}(M),[M]\rangle
+\sum^n_{k=1}\left(-\frac{1}{2}\right)^k \langle exp(w_k){\rm
Td}([c_k(\tau(M))]),[c_k(\tau(M))]\rangle.
$$
Therefore, the statement of the theorem is verified.


\section{Cancellation Formula
 of Alvarez-Gaume and Witten}

The gravitational anomaly cancellation formula~(\ref{MirCan}) of Alvarez-Gaume and Witten~\cite{AlWi} can also be obtain from Theorem~\ref{formula} by a suitable substitution. Let us consider a 12-dimensional compact stably almost complex manifold $M$ without boundary. Chose values of variables $z_i$, $i=1,2,\dots,$ in such a way that
$$
\prod^{\infty}_{k=1}(1-z_k)=\sqrt{1+\frac{1}{4}x^2},
$$
we also define variables $y_j$, $j=1,2,\dots,$ to satisfy
$$
\prod^{\infty}_{j=1}(1-y_jx)^{-1}=\frac{{\rm sinh}(x/2)}{x/2}.
$$
Using Theorem~\ref{formula} with this choice of $z_i$, $y_j$, $i,j=1,2,\dots,$ and
$t=0$, we obtain:
$$
\frac{1}{2^6}L(M)=\hat A(M)+\frac{1}{2^3}\hat A\left(\left[P_1(\tau(M))\right]\right)+
$$
\begin{equation}
\label{reduct}
+\frac{1}{2^7}\hat A\left(4\left[P_2(\tau(M)\right]-\left[P^2_1(\tau(M))\right]\right)
+\frac{1}{2^{10}}\left(p_{111}-4p_{12}+8p_3\right),
\end{equation}
where $\left[P_1(\tau(M))\right]$, $\left[P_2(\tau(M))\right]$, $\left[P^2_1(\tau(M))\right]$ are manifolds Poincare dual to the corresponding
tangent Pontryagin classes in cobordism, and
$p_{111}$, $p_{12}$, $p_3$ are the tangent Pontryagin numbers.

From the same substitution it is easy to see that
$$
\hat A\left(\left[P_1(\tau(M))\right]\right)=
\langle\frac{\partial}{\partial t}
\prod^6_{i=1}\frac{x_i\sqrt{1+t\left(e^{x_i}+e^{-x_i}-2\right)}}{{\rm sinh}(x_i/2)}
\left|_{t=0}\right.,[M]\rangle=
$$
$$
=
A(M,T)-12\hat A(M),
$$
where $x_i$, $i=1,\dots 6$, are the tangent Chern roots of the 12-dimensional 
stably almost complex manifold $M$, and
 $A(M,T)$ is the twisted $\hat A$-genus
defined by 
$$
\langle \prod_jA(x_j)\left(\sum_j
e^{x_j}+e^{-x_j}\right),[M]\rangle.
$$

Observe that if $W$ is a 4-dimensional stably almost complex manifold then
$\hat A(W)={\rm \bf s}_2([W])/24$, where ${\rm \bf s}_2$ is the second
Landweber-Novikov operation acting on $\Omega^*_U(\{pt\},\emptyset)$,
and $[W]\in \Omega^4_U(\{pt\},\emptyset)$ is the cobordism class defined by
$W$. Using standard rules of the Landweber-Novikov algebra, we deduce that
$$
{\rm \bf s}_2\left(4\left[P_2(\tau(M)\right]-\left[P^2_1(\tau(M))\right]\right)=
12p_{12}-24p_3-3p_{111}.
$$
Collecting similar terms in~(\ref{reduct}) we arrive at Alvarez-Gaume and Witten formula.
\begin{Cor}
For a 12-dimensional compact stably almost complex manifold without boundary
$$
L(M)=8A(M,TM)-32\hat A(M).
$$
\end{Cor}
\begin{Rem}
This formula is true for any compact oriented 12-dimensional manifold without boundary, as it only depends on Pontryagin classes.
\end{Rem}


\section{Toric Manifolds}

A toric variety is a closure of an orbit of an algebraic action of
$T^n$, and the quotient of the variety with respect to this action
is a simple polytope $P\subset {\mathbb R^n}$ (recall that a
polytope is simple if at every its vertex exactly $n$ faces of
codimension one meet each other). Toric varieties are not always
smooth. Polytopes for which the corresponding toric variety is
smooth, are regular (or Delzant) polytopes and varieties
themselves in this case are toric manifolds. More precisely, a
simple lattice polytope is regular if the edges emanating from
each vertex lie along vectors that generate the lattice $\mathbb
Z^n$. Each face $F_i$, $i=1,\dots, m$, of codimension one of a
regular polytope $P$ corresponds to a submanifold $M_{F_i}$ of
codimension two in the toric manifold $M_P$, $\dim M_P=2n$. The
torus action $T^n\times M_P\to M_P$ has a one-dimensional subgroup
$S^{F_i}\subset T^n$ which keeps $M_{F_i}$ fixed. $S^{F_i}$ can be
described by a primitive vector
$\lambda_i=(\lambda^1_i,\dots,\lambda^n_i)$ of the lattice
$\mathbb Z^n\subset \mathbb R^n$ which is inward normal to the
face $F_i\subset P\subset\mathbb R^n$. At each vertex $p\subset P$
we get $n$ vectors $\lambda_{i_j}$, $j=1,\dots,n$, corresponding
to each of the faces in the intersection:
$$
p=F_{i_1}\bigcap\cdots \bigcap F_{i_n}.
$$
Together, $\lambda_{i_j}$, $j=1,\dots, n$, form a basis of
$\mathbb Z^n\subset \mathbb R^n$.

There is also a smooth manifold ${\cal L}_P$ of dimension $n+m$
where $m$ is the number of codimension one faces of $P^n$ such
that $M^{2n}$ is a quotient of ${\cal L}_P$ over the action of a
torus $T^{m-n}$. The manifold ${\cal L}_P$ can be described as
follows. Let us list all codimension one faces of $P$:
$$
{\cal F}=\{F_1,F_2,...,F_m\}.
$$
For every face $F_i$ in ${\cal F}$ we define the corresponding
one-dimensional coordinate subgroup of $T^{\cal F}=T^m$ by
$T^{F_i}$. Now for a general face $G$ of $P$ let
$$
T^G=\prod_{G\subset F_i} T^{F_i}\subset T^{\cal F}.
$$
\begin{Def}
For a given simple polytope we define
$$
{\cal L}_P=\left(T^{\cal F}\times P^n\right)/\sim,
$$
where $(t_1,p)\sim (t_2,q)$ if and only if $p=q$ and
$t_1t^{-1}_2\in T^{G(q)}$.
\end{Def}
(In this definition $G(q)$ is the face of $P$ which contains $q$
strictly in its interior.)

For a toric manifold $M_P$ we define the characteristic
homomorphism $\lambda:T^{\cal F}\to T^n$, by means of isomorphisms
$$
\lambda (T^{F_i})\stackrel{\lambda_i}\cong S^{F_i}.
$$
This is a characteristic map in the sense of~\cite{BuPa,
DavisJan}. The kernel $H(\lambda)$ of $\lambda$ is
$(m-n)$-dimensional subgroup of $T^{\cal F}$ and as such it acts
on ${\cal L}_P$. This action is free and the factor space
$$
{\cal L}_P\stackrel{H(\lambda)}\longrightarrow {\cal
L}_P/H(\lambda)\cong M_P
$$
is the corresponding toric manifold. For each co-dimension one
face $F_i$ of the polytope $P$ we define a representation:
$$
\rho_i:T^{\cal F}\to GL(1,\mathbb C)\cong \mathbb C^*
$$
by means of a projection $T^{\cal F}\to T^{F_i}\subset
GL(1,\mathbb C)$. Consider a one-dimensional complex line bundle
$L_i$ over $BP=ET^{\cal F}\times_{T^{\cal F}} {\cal L}_P$ given by
$$
L_i=ET^m\times_{(T^{\cal F},\rho_i)}\left({\cal L}_P\times \mathbb
C\right).
$$
The following theorem~\cite{DavisJan} describes the equivariant
cohomology of ${\cal L}_P$.
\begin{Th} The space $ET^m\times_{T^{\cal F}}{\cal L}_P$ is
homotopy equivalent to $ET^n\times_{T^n} M^{2n}_P$. The cohomology
ring
$$
H^*(ET^m\times_{T^{\cal F}}{\cal L}_P,\mathbb Z)\cong
H^*(ET^n\times_{T^n} M^{2n}_P,\mathbb Z)
$$
is isomorphic to the face ring of $P$
$$
\mathbb Z[v_1,\dots, v_m]/I,\quad v_i=c_1(L_i).
$$
Under the projection $\pi_P:ET^n\times_{T^n} M^{2n}_P\to BT^n$ the
generator
$$
u_j\in H^*(BT^n,\mathbb Z)\cong \mathbb Z[u_1,\cdots, u_n]
$$
goes into
$$
\pi^*_P(u_j)=\lambda^j_1v_1+\cdots+\lambda^j_mv_m\in
H^*(ET^n\times_{T^n} M^{2n}_P,\mathbb Z).
$$
\end{Th}

\

At each vertex $p=F_{i_1}\bigcap\cdots\bigcap F_{i_n}$ of $P$ we
get $n$ linear independent vectors
$\lambda_{i_1},\dots,\lambda_{i_n}$ of $\mathbb Z^n$. Let
$M^t_p=\left(\mu_{p,i_1},\dots,\mu_{p,i_n}\right)$ be an $n\times
n$ column--matrix
($ \mu_{p,j}=\left(\mu^1_{p,j},\dots,\mu^n_{p,j}\right)^t$) such
that
$$
M_p\cdot \Lambda_p=I_n,\quad \Lambda_p=
\left(\lambda_{i_1},\dots,\lambda_{i_n}\right),
$$
and $I_n$ is the $n\times n$-identity matrix. Pick a generic
vector $\vec u\in\mathbb Z^n$ which is not orthogonal to any of
$\mu_{p,j}$.

\begin{Lemm}
\label{Localis} Under the Gysin map
$$
\left(\pi_P\right)_!: H^*(ET^n\times_{T^n} M^{2n}_P,\mathbb Z)\to
H^*(BT^n,\mathbb Z)
$$
the $n$th power of the first Chern class of $L_i$ goes into
\begin{equation}
\label{Gys1} \sum_{{i_1<\cdots<i_{n-1}\neq i}\atop {p_I=F_i\bigcap
F_{i_1}\bigcap\cdots \bigcap F_{i_{n-1}}\neq\emptyset}}
\frac{\langle \mu_{p_I,i}, \vec u \rangle^{n-1}}
{\prod^{n-1}_{j=1}\langle\mu_{p_I,i_j},\vec u\rangle },
\end{equation}
which is independent of the choice of a generic vector $\vec u$.
\end{Lemm}
{\bf Proof.} This is an elementary consequence of the fixed point
theory for compact Lie group actions~\cite{Audin,Pan,Quillen}. We
write down the fixed point formula for the action of $T^n$ on
$M_P$. At the fixed point corresponding to the vertex $p_I$ we get
a simple expression in terms of $v_i,v_{i_1},\dots, v_{i_{n-1}}$:
$$
\frac{v^n_i}{v_i\cdot v_{i_1}\cdots v_{i_{n-1}}}.
$$
Rewriting it uniformly in terms of $u_1,\cdots,u_n\in
H^*(BT^n,\mathbb Z)$ and substituting $u_j=\vec u_j$ leads to the
expression in the statement.

\

Let $\omega=(w_1\ge w_2\ge \cdots 0)$ be a partition of $n$ of
length $l$. Denote by $c_{\omega}(\eta)$ the product of Chern
classes $c_{w_1}(\eta), c_{w_2}(\eta),\dots $. Similar to
Lemma~\ref{Gys1} we deduce
\begin{Th}
For a smooth toric variety $M^{2n}$ its Chern number
$c_{\omega}(M)$ corresponding to a partition $\omega$ of $n$ can
be computed in terms of the fixed point data by means of a
formula:
$$
c_{\omega}(M)= \sum_{
 {
I_1=\{i_1< \cdots < i_{l}\}}\atop {I_2=\{i_{l+1}<\cdots <i_n\}
 }}
 \sum_{\sigma\in S_{l(\omega)} } \frac{
\prod^{l}_{j=1}\left(\mu^1_{p_I,i_j}u_1+\cdots+\mu^n_{p_I,i_j}u_n\right)^{w_{\sigma(j)}-1}}
{\prod^{n}_{j=l+1}\left(\mu^1_{p_I,i_j}u_1+\cdots+\mu^n_{p_I,i_j}u_n\right)}
$$
where the first summation is taken over $I_1 \bigcap
I_2=\emptyset$ and with
$$
p_I=F_{i_1}\bigcap \cdots \bigcap F_{i_{l}}\bigcap
F_{i_{l+1}}\bigcap \cdots \bigcap F_{i_n}\neq\emptyset.
$$
The formula is independent of the choice of a generic vector
$\vec u=(u_1,\dots,u_n)\in\mathbb R^n$.
\end{Th}

In dimension three all formulae can be rewritten just in terms of
$\lambda_j$ because
$$
\mu_{p,i}=\pm \frac{[\lambda_{p,j},\lambda_{p,k}]}{det\Lambda_p}.
$$
In particular,
\begin{equation}
\label{V3}
\left(\pi_P\right)_!(c_1(L_i))=\sum_{{i_1<i_2\neq
i}\atop {p_I=F_i\bigcap F_{i_1}\bigcap F_{i_2}\neq\emptyset}}
\frac{\langle \lambda_{p_I,i_1},\lambda_{p_I,i_2}, \vec u \rangle^2}
{\langle \lambda_{p_I,i},\lambda_{p_I,i_1}, \vec u \rangle\langle
\lambda_{p_I,i_2},\lambda_{p_I,i}, \vec u \rangle}
\end{equation}
 where
for three vectors $\vec x,\vec y,\vec z\in \mathbb R^3$ we denote the
corresponding oriented volume by $\langle \vec x,\vec y,\vec z\rangle$.


\section{Pick's Theorem}

In this section we apply Theorem~\ref{formula} for the case of
toric manifolds and its Chern classes. Recall that the $m$-th
tangent Chern class $c_m(\tau(M))$ of a toric manifold $M=M_P$
associated to a simple polytope $P$ of dimension $n$ is given
by a dual cohomology class to the sum of fundamental
classes of toric submanifolds $M_F$ associated to
$(n-m)$-dimensional faces $F=F^{(n-j)}_j$ of $P$. Let $w_P$ be the
Kahler form on $M_P$. It is well known~\cite{Guillemin} that the
equivariant class of the canonical Kahler form $w_P$ on $M_P$ is
$$
-a_1v_1-\cdots-a_m v_m\in H^2(ET^n\times_{T^n} M^{2n}_P,\mathbb
Z),
$$
where $a_j$, $j=1,...,m$, are defined by means of a system of
inequalities:
$$
x\in P\iff \langle x,\lambda_i\rangle \ge a_i,\quad i=1,\dots,m.
$$
We also set $x_j=i^*(v_j)$ where
$$
i:M^{2n}\subset ET^n\times_{T^n} M^{2n}\approx
ET^m\times_{T^m}{\cal L}_P
$$
is an embedding on any fibre.
\begin{Cor}
\label{Pick} The sum over all faces of a regular lattice polytope
$P\subset \mathbb R^n$ of the numbers of integral points inside
the closures of the faces taken with weights $(-1/2)^{n-\dim F}$,
is the twisted signature of the corresponding toric variety $M_P$:
\begin{equation}
\label{PiC} \langle exp(w_P)\prod^m_{i=1}
\frac{x_i/2}{\tanh(x_i/2)},[M_P]\rangle= \sharp(P\bigcap \mathbb
Z^n)+\sum^n_{k=1}\left(-\frac{1}{2}\right)^k \sum_{{F\subset
P}\atop {\dim F=(n-k)}}\sharp\left(F\bigcap \mathbb Z^n\right),
\end{equation}
where $x_i$ are the tangent Chern roots of $M_P$ ($m$ is the
number of codimension one faces).
\end{Cor}
{\bf Proof.} This is a 
straightforward consequence of
Theorem~\ref{Mirac} as due to~\cite{Danilov} the value
$$
\langle exp(w_F){\rm Td}(M_F),[M_F]\rangle
$$
is exactly the number of integral points in the face $F$.
\begin{Rem}
Formula~(\ref{PiC}) from Corollary~\ref{Pick} is valid for any simple lattice polytope. In order to verify this one needs to apply an equivariant 
generalization of
Theorem~\ref{Mirac} to the fibre bundle 
$ET^m\times_{T^m} {\cal L}_P\to BT^n$.
\end{Rem}
\begin{Exam}\rm{
Let $n=2$ and $P$ be a convex polygon in $\mathbb R^2$ with $m$
integral vertexes. Then the left hand side of (\ref{PiC}) is:
$$
\langle exp(w_P)\prod^m_{i=1}
\frac{x_i/2}{\tanh(x_i/2)},[M_P]\rangle= \langle
\frac{w^2_P}{2},[M_P]\rangle + \langle \prod^m_{i=1}
\frac{x_i/2}{\tanh(x_i/2)},[M_P]\rangle=
$$
$$
= Area(P)+\frac{\sigma(M_P)}{4},
$$
where $\sigma(M_P)$ is the signature of the corresponding toric
manifold. For any toric variety
$$
\sigma(M_P)=h_0(P)-h_1(P)+\cdots+(-1)^nh_n(P)=(-1)^nh_P(-1),
$$
where $(h_0,h_1,...,h_n)$ is the $h$-vector of $P$ and
$$
h_P(t)=h_0t^n+\cdots +h_n
$$
is the $h$-polynomial~\cite{Leu}. In particular, for $n=2$
$$
\sigma(M_P)=h_0(P)-h_1(P)+h_2(P)=4-m.
$$
Combining these observations with Corollary~\ref{Pick} we get
\begin{equation}
\label{Pick2}
Area(P)+\frac{4-m}{4}=(Int+Bd)-\frac{Bd+m}{2}+\frac{m}{4}=
Int+\frac{Bd}{2}-\frac{m}{4},
\end{equation}
where $Int$ is the number of integral points inside polygon $P$
and $Bd$ is the number of integral points on the boundary of $P$.
It is easy to see that (\ref{Pick2}) is equivalent to the
classical Pick's theorem:
$$
Area(P)=Int+\frac{Bd}{2}-1.
$$}
\end{Exam}
\begin{Exam}
\rm{Let $n=3$ and $P$ be a regular (Delzant) tetrahedron in $\mathbb R^3$
with $2$-dimensional faces $F_1,F_2,F_3,F_4$. Then the left hand
side of (\ref{PiC}) is
$$
\langle exp(w_P)\prod^4_{i=1}
\frac{x_i/2}{\tanh(x_i/2)},[M_P]\rangle= \langle
\frac{w^3_P}{6},[M_P]\rangle + \langle w_P\prod^4_{i=1}
\frac{x_i/2}{\tanh(x_i/2)},[M_P]\rangle=
$$
$$
= Vol(P)-\sum^4_{j=1} a_j\langle \frac{x_j}{3}\sum^4_{i=1}
\frac{x^2_i}{4},[M_P]\rangle= 
$$
$$
=
Vol(P)-\sum^4_{j=1}
\frac{a_j}{4}\langle \sum_{i\neq j}
\frac{x^2_i}{3},[M_{F_j}]\rangle-\sum^4_{j=1} a_j\langle
\frac{x^3_j}{12},[M_P]\rangle=
$$
$$
=Vol(P)-\sum^4_{j=1} \frac{a_j}{4}h_{F_j}(-1)- \sum^4_{j=1}
a_j\langle \frac{x^3_j}{12},[M_P]\rangle.
$$
Using~(\ref{V3}) we obtain for each fixed $j=1,\dots,4$:
$$
\langle x^3_j,[M_P]\rangle=\sum_{{i_1,i_2\neq j}\atop
{p=F_j\bigcap F_{i_1}\bigcap F_{i_2}\neq\emptyset}} \frac{\langle
\lambda_{p_I,i_1},\lambda_{p_I,i_2}, \vec u \rangle^2} {\langle
\lambda_{p_I,j},\lambda_{p_I,i_1}, \vec u \rangle\langle
\lambda_{p_I,i_2},\lambda_{p_I,j}, \vec u \rangle},
$$
where $\vec u\in\mathbb Z^3\in \mathbb R^3$ is an arbitrary vector. For
some fixed $i_1, i_2\neq j$ chose $\vec u=\lambda_{i_3}$, $i_3\neq
i_1,i_2,j$, because
$\lambda_{\alpha},\lambda_{\beta},\lambda_{\gamma}$ form a basis
of $\mathbb Z^3$ for any choice of pairwise different indices
$\alpha,\beta,\gamma$ we obtain that
$$
\langle x^3_j,[M_P]\rangle=\frac{\langle
\lambda_{p_I,i_1},\lambda_{p_I,i_2}, \lambda_{p_I,i_3} \rangle^2}
{\langle \lambda_{p_I,j},\lambda_{p_I,i_1}, \lambda_{p_I,i_3}
\rangle\langle \lambda_{p_I,i_2},\lambda_{p_I,j},
\lambda_{p_I,i_3} \rangle}=1.
$$
Thus,
$$
\langle exp(w_P)\prod^4_{i=1}
\frac{x_i/2}{\tanh(x_i/2)},[M_P]\rangle=Vol(P)-\sum^4_{j=1}\frac{a_j}{4}
-\sum^4_{j=1}\frac{a_j}{12}.
$$
 Let $Int$, $Fac$, $Edg$ be the numbers of integral points
inside polyhedron $P\subset \mathbb R^3$,  all faces of $P$ and
all edges of $P$ respectively. If $Vert$ is the number of vertexes
of $P$ then the right hand side of (\ref{PiC}) becomes
$$
Int+Fac+Edg+Vert-\frac{Fac+2Edg+3Vert}{2}+\frac{Edg+3Vert}{4}-\frac{Vert}{8}=
$$
$$
=Int+\frac{Fac}{2}+\frac{Edg}{4}+\frac{Vert}{8}.
$$
\begin{Cor}
\label{PiTe} For a regular (Delzant) tetrahedron in $\mathbb R^3$:
$$
Int+\frac{Fac}{2}+\frac{Edg}{4}+\frac{Vert}{8}=Vol(P)-\sum^4_{j=1}\frac{a_j}{3}.
$$
\end{Cor}
}
\end{Exam}

\section{Conclusion}
Applying the new cancellation formula for the twisted signature
and Todd genera of smooth toric varieties we derived a Pick's type
theorem for simple lattice polytopes. Explicit calculations by means of 
this formula are possible only in low-dimensional examples due to the lack of
a simple combinatorial discription for the $L$-class of toric varieties.
Finding such a description could lead to an interesting new invariant
of simple polytopes.

\

\noindent {\bf Acknowledgements.} The author expresses deep
gratitude to the anonymous referee for numerous valuable comments and
to Misha Feigin for a useful discussion.

\end{document}